\documentclass[12pt,a4paper]{amsart}
\theoremstyle{plain}
\usepackage{enumerate,amsfonts,amssymb,mathrsfs,amscd}
\usepackage{amsmath}
\usepackage[english]{babel}
\usepackage[utf8]{inputenc}
\usepackage{amsmath, amsthm, amssymb}
\usepackage[pdftex]{graphicx,color}

\usepackage{tikz,pgf}
\usetikzlibrary{calc}
\usepackage[colorlinks=true,linkcolor=red,citecolor=red]{hyperref}

\usepackage{etex}

\usepackage[shortlabels]{enumitem}
\usepackage{setspace}
\usepackage{tabularx,multirow}
\usepackage[arrow, matrix, curve]{xy}
\usepackage{MnSymbol}

\advance\hoffset-20mm \advance\textwidth40mm

\usepackage[normalem]{ulem}

\theoremstyle{plain}

\newtheorem{thm}{Theorem}[section]
\newtheorem{cor}[thm]{Corollary}
\newtheorem{lem}[thm]{Lemma}
\newtheorem{prop}[thm]{Proposition}

\theoremstyle{definition}
\newtheorem{defi}[thm]{Definition}

\newtheorem{conj}[thm]{Conjecture}

\newtheorem{rem}[thm]{Remark}
\newtheorem{rems}[thm]{Remarks}
\newtheorem{exa}[thm]{Example}

\newcounter{caseinproof}
0

\newcommand{\Der}{\operatorname{Der}}
\def\AA{{\mathbb A}}

\def\ZZ{{\mathbb Z}}
\def\QQ{{\mathbb Q}}

\newcommand{\A}{{\mathbb A}}

\def\kk{{\bf{k}}}
\def\GG{{\mathbb{G}}}
\def\CC{{\mathbb C}}

\def\PP{{\mathbb P}}

\def\cO{{\mathcal{O}}}

\newcommand{\p}{{\partial}}

\def\span{\mathop{\span}}
\def\Aut{\mathop{\rm Aut}}

\def\Aff{\mathop{\rm Aff}}
\def\SAut{\mathop{\rm SAut}}

\def\lie{\mathop{\rm lie}}

\def\SL{\mathop{\rm SL}}

\def\Bir{\mathop{\rm Bir}}

\def\Der{\mathop{\rm Der}}
\def\Jonq{\mathop{\rm Jonq}}
\def\jonq{\mathop{\rm jonq}}

\def\deg{\mathop{\rm deg}}

\def\exp{\mathop{\rm exp}}

\def\span{\mathop{\rm span}}

\def\ML{\mathop{\rm ML}}
\def\GL{\mathop{\rm GL}}

\def\lie{\mathop{\rm lie}}

\def\ll1{l_{\lambda}^{-1}(1)}
\def\lm1{l_{\mu}^{-1}(1)}

\newcommand\qmatrix[2][1]{\left(\renewcommand\arraystretch{#1}
\begin{equation}gin{array}{*{20}r}#2\end{array}\right)}

\begin{document}

\title{Tits-type alternative for certain groups acting on algebraic surfaces}
\author{Ivan~Arzhantsev and  Mikhail~Zaidenberg}
\address{National Research University Higher School of Economics, Faculty of Computer Science,
Pokrovsky Boulevard 11, Moscow, 109028 Russia}
\email{arjantsev@hse.ru}

\address{Institut Fourier, UMR 5582, Laboratoire de Math\'ematiques\\
Universit\'e Grenoble Alpes, CS 40700, 38058 Grenoble cedex 9, France
}
\email{Mikhail.Zaidenberg@univ-grenoble-alpes.fr}
\date{}

\thanks{The first author was supported by the grant RSF-DST 22-41-02019}.
\begin{abstract} 
According to a theorem of Cantat and Urech, an analog of the classical Tits alternative holds for the group of birational automorphisms of a compact complex K\"{a}hler surface.
We established in~\cite{AZ21} the following Tits-type alternative: if $X$ is a toric affine variety and $G\subset\Aut(X)$  is a subgroup  generated by a finite set 
of unipotent subgroups normalized by the acting torus then either $G$ contains a nonabelian free subgroup or $G$ is a unipotent affine algebraic group. 
In the present paper we extend the latter result to 
any group $G$ of automorphisms of a complex affine surface
generated by a finite collection of unipotent algebraic subgroups. It occurs that either 
$G$ contains a nonabelian free subgroup or $G$ is a metabelian unipotent algebraic group.

\end{abstract}
\maketitle

\thanks{
{\renewcommand{\thefootnote}{} \footnotetext{2020
\textit{Mathematics Subject Classification:}
14J50, 14R20, 14L30, 14E07,\,22F50.\mbox{\hspace{11pt}}\\{\it Key words}:  affine surface, group action, unipotent group, one-parameter subgroup, Cremona group, Tits alternative.}}

{\footnotesize \tableofcontents}

\section*{Introduction} 
Let $\kk$ be an algebraically closed field of characteristic zero, and let 
$\GG_{\mathrm{a}}$ and  $\GG_{\mathrm{m}}$ stand for the additive and multiplicative groups of $\kk$, respectively, viewed as algebraic groups. 
Given an algebraic variety $X$ over $\kk$ and a subgroup $U$ of $\Aut (X)$ isomorphic to $\GG_{\mathrm{a}}$ 
and acting regularly on $X$, $U$ is called a \emph{one-parameter unipotent subgroup} or a \emph{$\GG_{\mathrm{a}}$-subgroup}.  
Any $\GG_{\mathrm{a}}$-subgroup has the form $U=\exp(t\p)\subset\SAut(X)$ where $t\in\kk$ and $\p$ 
is a locally nilpotent derivation (an \emph{LND} for short) of the structure ring $\cO(X)$.

The general linear group $\GL_n(\kk)$ verifies the Tits alternative~\cite{Tit72}, that is, any subgroup $G\subset\GL_n(\kk)$
either contains a nonabelian free subgroup or is a finite extension of a solvable group 
\footnote{In~\cite{AZ21} we called this property the \emph{enhanced Tits alternative} in order to distinguish 
from the usage of the term 
``Tits alternative'' applied exceptionally to finitely generated subgroups.}. 
 In the case of a connected linear algebraic group $G$ defined over $\kk$ this alternative follows from the 
 Levi-Malcev theorem. Indeed, consider
the Levi decomposition $G=R_{\mathrm u}\rtimes L=R_{\mathrm u}\rtimes (S\cdot T)$ where $R_{\mathrm u}$ 
is the unipotent radical of $G$, $L$ is a reductive Levi subgroup, 
$S$ is a semisimple algebraic group and $T$ is an algebraic torus. We have the alternative:  
either $G$ is solvable or $G$ contains a semisimple algebraic subgroup $S$ of positive dimension. In the  latter case
$S$  contains a subgroup isomorphic  either to ${\rm SL}(2,\kk)$ or to ${\rm PSL}(2,\kk)$, which in turn contains a nonabelian free subgroup.
It follows that $G$ has no nonabelian free subgroup if and only if $S$ is trivial, if and only if $G$ is solvable. 

If $G$ has no nontrivial character then $L=S$, that is, $G=R_{\mathrm u}\rtimes S$. In this case $G$ is generated by   
a  finite collection of $\GG_{\mathrm{a}}$-subgroups normalized by a given  maximal torus $T_{\rm max}$ of $S$. 
Indeed, this clearly holds for  $R_{\mathrm u}$. In turn, $S$ is generated by  the finite collection of root $\GG_{\mathrm{a}}$-subgroups with respect to $T_{\rm max}$. 
 Such a group $G$ has no nonabelian free subgroup if and only if $S$ is trivial, if and only if $G$ is unipotent.
The latter issue holds as well in the following more general setup.

\begin{thm}[{\rm \cite[Theorem~1.1]{AZ21}}]\label{thm:0.4} Consider a toric affine variety $X$ over $\kk$. 
\footnote{In~\cite{AZ21} $X$ is assumed to be a toric variety with no toric factor. However, the statement obviously holds without the latter assumption.}
Let a subgroup  $G=\langle U_1,\ldots,U_s\rangle$ of $\Aut(X)$ 
be generated by a finite collection of $\GG_{\mathrm{a}}$-subgroups $U_1,\ldots,U_k$ normalized by the acting torus. 
Then  either $G$ contains a nonabelian free subgroup or $G$ is a unipotent algebraic group.
\end{thm}

Assuming that $\dim X\ge 2$, $X$ is smooth in codimension 2 and has no nonconstant invertible regular function one can always find a group $G$ 
as in Theorem~\ref{thm:0.4} acting on the smooth locus of $X$ highly transitively, that is, $n$-transitively for any natural number $n$~\cite[Theorem~1.3]{AKZ19}. 
If $G$ acts doubly transitively then $G$ contains a nonabelian free subgroup, see~\cite[Corollary~1.2]{AZ21}. 

The natural question arises whether the assertion of Theorem~\ref{thm:0.4} remains true for more general affine varieties,
disregarding the existence of a torus action. 
In Section~\ref{sec:aff-surf} 
we give a positive answer 
for the complex affine surfaces. 
Our main result (Theorem~\ref{thm:aff-surf}) is as follows. 

\begin{thm}\label{thm:main}
Let $X$ be an affine algebraic surface over $\CC$ and
$G$ be a subgroup of 
$\Aut(X)$ generated by a finite collection of  $\GG_{\mathrm{a}}$-subgroups $U_1,\ldots,U_k$. 
Assume that $G$ contains no nonabelian free subgroup. Then $G$ is a metabelian unipotent affine algebraic group. 
\end{thm}

Recall that any abelian group is also metabelian. 
The proof of  Theorem~\ref{thm:main} uses Propositions~\ref{thm:A2} and~\ref{thm:proj-surf}. 
The former claims the validity of the theorem for the affine plane $X=\A^2_\CC$. Notice that 
the conclusion 
of Theorem~\ref{thm:main} does not hold any longer for the Cremona group $\Bir(\PP^2_\kk)$, see Example~\ref{exa:inf-dim}. 
In Proposition~\ref{thm:proj-surf} we describe the exceptional cases. The proof of Theorem~\ref{thm:main} is done in Section~\ref{sec:aff-surf}. In Section~\ref{sec:unip-grps}
we study the actions of unipotent algebraic groups on affine surfaces. 
Using Theorem~\ref{thm:main} we deduce a Tits-type alternative for Gizatullin affine surfaces, see  Propositions~\ref{prop:unip}  and \ref{new-cor}. 

The reason why a solvable subgroup $G$  of $\Aut(X)$ generated by a finite set of unipotent algebraic subgroups is unipotent is clear in the case where  $G$
is contained in an affine algebraic group $\tilde G$. Indeed, 
in this case $G$ is contained in a Borel subgroup $B$ of $\tilde G$ and  any unipotent subgroup of $G$ is contained in the unipotent radical of $B$. Hence $G$  itself is unipotent. 

However, a priori it is not clear why such a subgoup $G$ must be  algebraic. The latter holds e.g. if $G$ consists of triangular automorphisms, see~\cite[Proposition 15.2.5]{FK18}, cf also \cite{AS20}.
This also happens in the setting of Theorem~\ref{thm:main}.
After our paper had been submitted, this was proven to be the case for any affine variety defined over $\kk$, see \cite[Theorem B]{KZ22}. 

The proofs in the present paper rely on several results established for varieties defined over $\CC$,  see \cite{Can11, Des07, Des15, Des21, Lam01, Ure21} and Section \ref{sec:prelim} for an account of the results from these papers used in the proofs.
We believe however that the  following analog of Theorem~\ref{thm:main} holds. 

\medskip

\begin{conj}\label{conj} \emph{Let $X$ be an affine algebraic variety defined over an algebraically closed field 
and let $G$ be
 a subgroup  of $\Aut(X)$ generated by a finite collection of  $\GG_{\mathrm{a}}$-subgroups $U_1,\ldots,U_k$. Then
either $G$ contains a nonabelian free subgroup, or $G$ is a unipotent affine algebraic group. }
\end{conj} 

\medskip

In Section \ref{sec:note} (added in proof) we provide 
an alternative proof of Theorem \ref{thm:main} valid for any algebraically closed field of characteristic zero. 
This proof applies a recent result from \cite{KZ22} and relies on the Tits alternative for the Cremona group of rank $2$ over $\kk$. 
The proof of the latter is 
sketched in the book [Lam22] under construction. 
Therefore, we consider for the moment our alternative proof to be conventional. 

\medskip

 {\bf Acknowledgments.} We are grateful to Yuri Prokhorov, Christian Urech and St\'ephane 
 Lamy for useful remarks and indications in the literature. Our thanks also go to 
 the anonymous reviewer for several suggestions that have improved the article.
 
\medskip

\section{Preliminaries}\label{sec:prelim}

\subsection{Subgroups of $\Aut(\A^2_\kk)$ and the Tits alternative}\label{ss:A2}
Recall the following facts.

\begin{thm}[Rentchler~\cite{Ren68}] \label{thm:R}
Any $\GG_{\mathrm{a}}$-subgroup of $\Aut(\A^2_\kk)$ is conjugate to a subgroup of triangular transformations
\[(x,y)\mapsto (x+\tau p(y), y), \quad \tau \in\kk\] 
generated by a triangular LND $\p=p(y)\p/\p x$ for some $p\in\kk[y]$.
\end{thm}

For the following theorem see e.g. 
\cite{Jun42, vdK53, Nag72, Kam75, Wri75, Wri78}. 

\begin{thm}[Jung-van der Kulk] \label{thm:JK}
There is a free amalgamated product decomposition $\Aut(\A^2_\kk)=A*_C B$ 
where $C=A\cap B$, $A=\Aff(\A^2_\kk)$ is the affine group of the affine plane and $\Jonq(\A^2_\kk)$ 
is the group of de Jonqui\`eres transformations
\[(x,y)\mapsto (\alpha x + p(y), \beta y + \gamma)\quad\text{with}\quad\alpha,\beta\in\kk^*,\,\,\,\gamma\in\kk.\]
\end{thm}

Notice that Theorem~\ref{thm:R} can be deduced from Theorem~\ref{thm:JK}. Indeed, 
any algebraic subgroup of $\Aut(\A^2_\kk)$ has bounded degree, hence also bounded length~\cite[1.10]{Wri79}. 
However, 
any  subgroup of bounded length of an amalgamated product $A*_C B$ is conjugate 
to a subgroup of one of the factors $A$ and $B$~\cite[Theorem~8]{Ser80}; 
cf.~\cite[Corollary~4.2 and Theorem~4.3]{Kam79} and~\cite[Propositions~0.35 and ~1.11]{Wri79} 
for alternative arguments. 

The abelian subgroups of an amalgamated product of two groups are described as follows. 

\begin{thm}[Moldavanskii~\cite{Mol67}]\label{thm:M} 
Consider a free product with amalgamated subgroup $F=A *_C B$. 
Then for any abelian subgroup $G$ of $F$ precisely one of the following holds:
\begin{enumerate}
\item[(1)] $G$ is conjugate in $F$ to a subgroup of either $A$ or $B$;
\item[(2)] $G =\bigcup_i H_i$ where
$H_0\subset H_1 \subset\ldots\subset H_i \ldots$ is an infinite non-stationary chain of subgroups 
such that each $H_i$ is conjugate 
in $F$ to a subgroup of $C$ while $G$ is not conjugate to any subgroup of $A$ or $B$;
\item[(3)] $G = H\times \langle g\rangle$ where $H$ is conjugate to a subgroup of $C$  
and $g\in G$ is an element of infinite order non-conjugate to any element of $A$ or $B$. 
\end{enumerate}
\end{thm}

See~\cite{Wri79} for an alternative approach based on the Bass-Serre theory
~\cite{Bas76, Ser80}. 
This approach is also used in~\cite{Lam01} and \cite{Lam22} to obtain the following result.
\begin{thm}[{\rm \cite{Lam01}, \cite[Proposition 16.14]{Lam22}}]\label{cor:L} 
The group $\Aut(\A^2_\kk)$ verifies the Tits alternative. 
\end{thm}
For $\kk=\CC$} a more detailed analysis is available, see~\cite[Theorem~2.4,~Propositions~3.12~and~4.10]{Lam01}. 
In fact, this can be extended to general fields, see Remark \ref{rem:Lamy-book}.
\begin{thm}[Lamy] \label{thm:L} For any subgroup 
$G$ of $\Aut(\A^2_\CC)$  one of the following holds: 
\begin{enumerate}
\item[(1)] either $G$ is conjugate to a subgroup of $\Aff(\A^2_\CC)$ or $\Jonq(\A^2_\CC)$,\\
\noindent or $G$ belongs to one of the following types:
\item[(2)] $G$ is an abelian torsion group;
\item[(3)] $G$ is a solvable group which contains a subgroup of finite index isomorphic to $\ZZ$;
\item[(4)] $G$ contains a nonabelian free subgroup.
\end{enumerate}
\end{thm}

Let us remind that the de Jonqui\`eres group $\Jonq(\A^2_\kk)$ is solvable, see e.g.~\cite[Corollary~2.5]{Lam01}. 
Notice that an abelian subgroup $G\subset\Aut(\A^2_\kk)$ of type (2) resp. (3) above fits in case (2) 
resp. (3) of Theorem~\ref{thm:M}. 

\subsection{Subgroups of the Cremona group $\Bir(\PP^2_\CC)$}\label{ss:Bir2}
More generally, we have the following theorem.

\begin{thm}[Cantat-Urech~\cite{Can11, Ure21}] \label{thm:CU}
Let $X$ be a compact complex  K\"{a}hlerian surface. Then the group $\Bir(X)$ verifies the Tits alternative. 
\end{thm}

The birational actions of algebraic groups were considered e.g. by Weil~\cite{Wei55}, 
Rosenlicht~\cite{Ros56}, Matsumura~\cite{Mats63}, Grothendieck~\cite[\S20.5]{Gro67} and Demazure~\cite[\S1]{Dem70}. 
 We address~\cite{Ser10},~\cite{BF13} and~\cite{Bla17} for a recent treatment. Recall the following definitions.

\begin{defi}\label{def:Bir} 
Let $X, Z$ be algebraic varieties over $\kk$.
A mapping $f\colon Z\to\Bir(X)$ is said to be a \emph{morphism} if the following holds:
\begin{itemize}
\item the induced map 
\[Z\times X\to Z\times X,\quad (z,x)\mapsto (z, f(z)(x))\]
is a birational transformation which induces a biregular isomorphism $U\stackrel{\cong}{\longrightarrow} V$ 
between some dense open subsets 
$U, V\subset Z\times X$ whose projections to $Z$ are surjective.
\end{itemize}
\end{defi}

\begin{defi}\label{def:act} 
Let $G$ be an algebraic group over $\kk$. A \emph{birational action} of $G$ on an algebraic variety $X$ over $\kk$ is a morphism 
$G \to \Bir(X)$ which is also a group homomorphism. 
A subgroup $H$ of $\Bir(X)$ is called \emph{algebraic} if $H$ is the image of some algebraic group $G$ 
under a homomorphism $G \to \Bir(X)$  which is a morphism. If $H\neq 1$ is algebraic and $G=\GG_{\mathrm{a}}$ 
then $H$ is called a $\GG_{\mathrm{a}}$-\emph{subgroup} of $\Bir(X)$.
\end{defi}

Recall that any algebraic subgroup of $\Bir(\PP^n_\kk)$ is a linear group~\cite[Remark~2.21]{BF13}.
The next lemma follows directly from the definition of algebraic
group action. 

\begin{lem}[{\rm \cite[\S 1]{Dem70}; see also~\cite[Remark~2.4]{Bla17}}] Let $X$ and $Y$ be algebraic varieties. 
Then any birational map $X\dasharrow Y$ induces a group isomorphism $\Bir(X)\cong\Bir(Y)$ 
and conjugates birational group actions on $X$ to birational group actions on $Y$. 
\end{lem}

For the next theorem see~\cite[Theorem on p.~375]{Wei55},~\cite[Theorem~1]{Ros56}, see also~\cite[Theorem~4.6]{Zai95}, \cite[Theorem~2.17]{Bla17}, \cite[Theorem~1]{Kra18} and~\cite{Bri22}.

\begin{thm}\label{thm:Weil} Let $X$ be an algebraic variety, $G$ be an algebraic group and $G \to \Bir(X)$ be a birational group action. 
Then, there exists a birational map $X \dasharrow Y$ where $Y$ is an algebraic variety that conjugates this action 
to a biregular group action of $G$ on $Y$.
\end{thm}

Next we turn to classification of nilpotent and solvable subgroups of the complex Cremona group $\Bir(\PP^2_\CC)$. 

\begin{thm}[D\'eserti~\cite{Des07}]\label{thm:D-nilp}
Let $G$ be a nilpotent subgroup of $\Bir(\PP^2_\CC)$. Then $G$ is either a torsion group
or a virtually metabelian group.
\end{thm}

The next theorem is a reformulation of~\cite[Theorem~8.49]{Des21}, 
which in turn integrates the results of~\cite{Des15} and~\cite[Theorem~7.1]{Ure21}
based in particular on~\cite[Proposition~6.14 
and Theorem~7.7]{Can11}.  

\begin{thm}[D\'eserti--Urech]\label{thm:DU-solv}
Let $G$ be a solvable subgroup of $\Bir(\PP^2_\CC)$. Then one of the following
holds:
\begin{enumerate}
\item[(1)]  
$G$ is conjugate to a subgroup of ${\rm PGL}(2, \CC(t)) \rtimes {\rm PGL}(2, \CC)$;
\item[(2)]  $G$ is conjugate to a subgroup of an affine algebraic group defined over $\CC$;
\item[(3)]  $G$ is conjugate to a subgroup of the normalizer $N(T)\cong\GG_{\mathrm m}^2\rtimes\GL(2, \ZZ)$  
of the standard $2$-torus $T\subset\Bir(\PP^2_\CC)$; 
\item[(4)]  Up to passing possibly to a subgroup of index $2$, $G$  
is isomorphic to $H\rtimes \ZZ$ for a finite subgroup $H$ of $\Bir(\PP^2_\CC)$. 
\end{enumerate}
\end{thm}

\begin{rem}\label{rem:der-length} By Theorem~\ref{thm:D-nilp} the derived length of 
any connected nilpotent subgroup of $\Bir(\PP^2_\CC)$ is at most 2.
 The maximal derived length of a solvable subgroup $G\subset\Aut(\A^2_\CC)$ 
 ($G\subset\Bir(\PP^2_\CC)$, respectively) equals 5~\cite[Proposition~3.14]{FP18}  
 (is bounded above by 8, respectively~\cite[Theorem~1.10]{Ure21}). 
 The derived length of a solvable closed connected subgroup of $\Bir(\PP^2_\CC)$ is at most $4$ and this bound is reached for the Borel subgroup 
 $\mathcal{B}_2 \cong {\rm Aff}_1(\CC(y)) \rtimes {\rm Aff}_1(\CC)$ \cite[Proposition 4 and Lemma 52]{FH21}.
See also ~\cite[Proposition~2.2]{Cor13} for an example of a solvable $3$-generated subgroup 
$G\subset\Bir(\PP^2_\CC)$ of derived length 3 which admits no faithful linear representation. 

For any $n\ge 1$ there is a nilpotent  2-generated torsion-free subgroup $G_n$ of 
$\Bir(\PP^2_\kk)$ of nilpotency class $n$~\cite[Lemma~2.1]{Cor13}.  Furthermore, the 
unipotent algebraic subgroup $G$ of $\Aut(\A^2_\CC)$ generated by the $\mathbb{G}_{\mathrm{a}}$-subgroups
\[U_1=\{(x,y)\mapsto (x+s, y)\}\quad\text{and}\quad U_2=\{(x,y)\mapsto (x, y+\tau x^n)\}, \quad s,\tau \in\CC\] has nilpotency class $n$.
\end{rem}

\section{The case of the affine plane}\label{sec:aff-plane}
In this section we prove Theorem~\ref{thm:main} in the particular case 
$X=\A^2_\CC$.  In section~\ref{sec:aff-surf} we reduce the general case to this one. 

\begin{prop}\label{thm:A2} 
Let 
$G$ be a subgroup of 
$\Aut(\A^2_{\CC})$ generated by a finite collection of 
$\GG_{\mathrm{a}}$-subgroups $U_1,\ldots,U_k$ defined over $\CC$. 
Assume that $G$ contains no nonabelian free subgroup. 
Then $G$ is a unipotent affine algebraic group conjugate in $\Aut(\A^2_\CC)$ to a subgroup of the metabelian
unipotent de Jonqui\`eres group 
\[{\jonq}(\A^2_{\CC})=\{(x,y)\mapsto (x+p(y), y+\gamma)\,|\, p\in\CC[y], \,\,\,\gamma\in\CC\}.\]  
\end{prop}

\begin{proof} Under our assumption $G$ is solvable, see Corollary~\ref{cor:L}. So, it is of one of types (1)-(3) in Theorem~\ref{thm:L}.
Since $G$ is a connected subgroup of the ind-group $\Aut(\A^2_{\CC})$ it cannot be of type (2) or (3). Hence, $G$ is of type (1),
that is, $G$ is conjugate to a subgroup of one of the groups $\Aff(\A^2_{\CC})$ and $\Jonq(\A^2_{\CC})$.

Let $G$ be conjugate to a subgroup of $\Aff(\A^2_{\CC})$. Then $G$ is a connected solvable algebraic group.  
Hence $G$ is contained in a Borel subgroup $B$ of $\Aff(\A^2_{\CC})$. The unipotent subgroups $U_1,\ldots,U_k$ 
are contained in the unipotent radical $U$ of $B$, hence also $G\subset U$. 
Consider the standard faithful representation $\Aff(\A^2_{\CC})\to \GL(3,\CC)$. 
Since the Borel subgroups of $\Aff(\A^2_{\CC})$ are conjugated we may suppose that the image of $G$ in $ \GL(3,\CC)$ 
consists of upper triangular matrices of order 3 with unit diagonal. 
In other words, $G$ is an algebraic subgroup of the metabelian unipotent affine algebraic group
\[(x,y)\mapsto (x+ay+b, y+\gamma)\quad\text{where}\quad a,b,\gamma\in\CC.\]
In this case $G$ is a metabelian unipotent algebraic group (of dimension at most 3).

Let now $G$ be conjugate to a subgroup of $\Jonq(\A^2_{\CC})$. We may assume that 
$G\subset\Jonq(\A^2_{\CC})$. There is a decomposition into a semidirect  product
\[{\Jonq}(\A^2_{\CC})={\jonq}(\A^2_{\CC})\rtimes T\]
 where $T=\{(x,y)\mapsto (\alpha x, \beta y)\,|\,\alpha, \beta\in\CC^*\}$  is the maximal torus. 
Since any homomorphism $U_i\to T$ is trivial $U_i$ is contained in  the unipotent radical
${\jonq}(\A^2_{\CC})$ of ${\Jonq}(\A^2_{\CC})$ for $i=1,\ldots,k$. 
Hence also $G\subset {\jonq}(\A^2_{\CC})$. There is the ascending filtration 
\[{\jonq}(\A^2_{\CC})=\bigcup_{d=1}^\infty {\jonq}_d(\A^2_{\CC})\quad\text{where}\quad 
 {\jonq}_d(\A^2_{\CC})=\{f\in {\jonq}(\A^2_{\CC})\,|\deg(f)\le d\}\] 
is a metabelian unipotent  affine algebraic group. This defines on $ {\jonq}(\A^2_{\CC})$ 
a structure of a nested ind-group. So $U_i\subset  {\jonq}_d(\A^2_{\CC})$ for all $i=1,\ldots,k$ and for some 
$d\gg 1$. Therefore, $G\subset  {\jonq}_d(\A^2_{\CC})$ is a metabelian unipotent affine algebraic group.
\end{proof}

\begin{rem}
The validity of the Tits alternative for the group $\Aut(\A^3_\CC)$ is still an open question. It holds, however, in the tame subgroup of $\Aut(\A^3_\CC)$, see \cite{LP22}.
See also ~\cite[Corollary~2]{Pop17} (cf. also~\cite[Corollary~3.2]{DF91}) for the following partial result.

\smallskip

\noindent \emph{Every connected solvable affine algebraic subgroup $G$ of $\Aut(\A^3_\kk)$ is rationally triangulable, 
that is,  conjugate in $\Bir(\A^3_\kk)$ to a subgroup of the de Jonqu\`eres group  ${\rm Jonq}(\A^3_\kk)$ of triangular automorphisms of $\A^3_\kk$.}

\smallskip

\noindent It follows that the derived length of $G$ is at most 4. Indeed, the derived length  of ${\rm Jonq}(\A^3_\kk)$ equals 4~\cite[Lemma~3.2]{FP18}. 
\end{rem}

\section{The case of the projective plane}
In this section we prove the following proposition.
\begin{prop}\label{thm:proj-surf} 
Let $X$ be an algebraic surface over $\CC$ and let 
$G$ be a subgroup of 
$\Bir(X)$ generated by a finite collection of $\GG_{\mathrm{a}}$-subgroups $U_1,\ldots,U_k$. 
Assume that $G$ contains no nonabelian free subgroup. Then $G$ is a solvable group. Furthermore, one of the following holds:
\begin{enumerate}
\item[(1)] $G$ is a metabelian unipotent affine algebraic group;
\item[(2)] $X$ is rational, $G$ preserves a pencil of rational curves 
on $X$ and contains a normal subgroup $N$ such that $G/N\cong\GG_{\mathrm{a}}$. 
\end{enumerate}
\end{prop}
Recall~\cite[Theorem~1.1]{Des21} that, given a pencil $\Lambda$ of rational curves on $\PP^2_\kk$,  the subgroup 
of the Cremona group $\Bir(\PP^2_\kk)$ consisting of transformations preserving $\Lambda$ is conjugate to the 
de Jonqui\`eres subgroup 
$J$ of the Cremona group $\Bir(\PP^1_\kk\times\PP^1_\kk)\cong\Bir(\PP^2_\kk)$. One has
\[J={\rm PGL}(2, \kk(t)) \rtimes {\rm PGL}(2, \kk)\] 
 where ${\rm PGL}(2, \kk)$ acts on $\PP^1_\kk\times\PP^1_\kk$ via automorphisms of the second factor 
 and acts on ${\rm PGL}(2, \kk(t))$ via the natural ${\rm PGL}(2, \kk)$-action on the rational function field $\kk(t)$. 

The proof of Proposition~\ref{thm:proj-surf} is preceded by the following lemmas.
\begin{lem}\label{lem:P1xP1} 
Let $G$ be a solvable subgroup of the subgroup ${\rm PGL}(2, \kk(t))\subset
J$ generated by a finite collection of $\GG_{\mathrm{a}}$-subgroups $U_1,\ldots,U_k\subset \Bir(\PP^1_\kk\times\PP^1_\kk)$ 
defined over $\kk$.
Then $G$ is an abelian subgroup of the unipotent radical $R_{\mathrm{u}}(\tilde B)$ where $\tilde B$ is a Borel subgroup  
of ${\rm PGL}(2, \kk(t))$. 
\end{lem}
\begin{proof} 
Let $\bar G$ be the minimal algebraic subgroup of ${\rm PGL}(2, \kk(t))$ defined over $\kk(t)$ which contains $G$. 
It is well known that $\bar G$ is solvable;  
see e.g. \cite[Lemma~2.3]{FP18}. 
It is also connected. Indeed, the component group of $\bar G$ is finite.
Since the $\GG_{\mathrm{a}}$-subgroups $U_i$ are infinitely divisible they are contained in 
the neutral component $\bar G^0$ of $\bar G$. 
 Hence $\bar G=\bar G^0$ is connected. 

It follows that $\bar G$ is contained in a Borel subgroup $\tilde B$ of ${\rm PGL}(2, \kk(t))$. 
Since every unipotent subgroup $U_i$ of $\tilde B$  is contained in the unipotent radical 
$R_{\mathrm{u}}(\tilde B)\cong \GG_{\mathrm{a}}(\kk(t))$,  then  so is $G$. 
\end{proof}
\begin{lem}\label{lem:alg} 
Let $X$ be a  projective variety over $\kk$ and $G$ be a solvable subgroup of 
$\Bir(X)$ generated by a finite collection of $\GG_{\mathrm{a}}$-subgroups $U_1,\ldots,U_k\subset \Bir(X)$. 
Assume that $G$ is conjugate to a subgroup of $\Aut(S)$ where $S$ is a projective variety. 
Then $G$ is a unipotent affine algebraic group. The same conclusion holds if $G$ is conjugate to a subgroup of an affine algebraic group.
\end{lem}
\begin{proof} Assume $G$ is conjugate to a subgroup of $\Aut(S)$. The group $G$ being connected the image of $G$ in $\Aut(S)$ is contained in the neutral component $\Aut^0(S)$. 
By the Matsusaka theorem~\cite{Mat58} $\Aut^0(S)$ is an algebraic group. 
According to a Chevalley
theorem~(see \cite[Theorem~16]{Ros56})
there exists a normal connected affine algebraic subgroup $H$ of $\Aut^0(S)$ such that $A=\Aut^0(S)/H$ is an abelian variety. 
Since any homomorphism $U_i\to A$ is trivial the image of $U_i$ in $\Aut^0(S)$ is contained in $H$ for $i=1,\ldots,k$. Hence also 
the image of $G$ is contained in $H$. 
Being solvable and connected, the image of $G$ is contained in a Borel subgroup $B$
of $H$. There is a decomposition  $B=U\rtimes T$ where $U$ is the unipotent radical of $B$ and $T$ is a maximal torus of $H$. 
Since any morphism $U_i\to T$ is trivial the image of $U_i$ is contained in $U$  for $i=1,\ldots,k$. Therefore the image 
of $G$ is a closed subgroup of $U$, and so, a unipotent affine algebraic group. The same argument proves the last assertion of the lemma.
\end{proof}
\begin{rem} Let $X$ be a normal affine variety over $\kk$ and let $H_1,\ldots,H_k$ be algebraic subgroups of $\Aut(X)$. 
If these subgroups are contained in an algebraic subgroup $H$  of $\Aut(X)$ then by the Sumihiro theorem~\cite[Theorem~3]{Sum74} 
there exists an $H$-equivariant completion $\bar X$ of $X$ so that the action of $H_1,\ldots,H_k$ on $X$ extends to their action on $\bar X$. 
The proof of Lemma~\ref{lem:alg} shows the partial converse, namely,

\smallskip

\noindent  \emph{Assume that $H_1,\ldots,H_k$ are connected and their action 
on $X$ can be extended to a completion $\bar X$ of $X$.
Then $G=\langle H_1,\ldots,H_k\rangle$ is a connected linear algebraic group.}

\smallskip

Recall that $G\subset\Aut^0(\bar X)$ is closed, hence an algebraic subgroup, see~\cite[Proposition~7.5]{Hum75}.
\end{rem}
\begin{lem}\label{lem:two actions} 
Let $C$ be a curve over $\kk$ and $U$ and $U'$ be two different $\GG_{\mathrm{a}}$-subgroups of $\Aut(C)$.
Then $C\cong\PP^1_\kk$, the fixed points of $U$ and $U'$ are distinct and $\langle U,\,U'\rangle\cong {\rm PGL}(2,\kk)$. 
\end{lem}
The proof is a simple exercise and will be omitted. 
\begin{lem}\label{lem:irrat} 
Let $X$ be a projective surface over $\kk$ and $G$ be a
subgroup of 
$\Bir(X)$ generated by a finite collection of $\GG_{\mathrm{a}}$-subgroups $U_1,\ldots,U_k$. 
If $X$ is irrational then either one of the subgroups $\langle U_i,\,U_j\rangle$  surjects onto 
${\rm PGL}(2,\kk)$ and so contains a nonabelian free subgroup,
or $G$ is an abelian unipotent algebraic group. The latter holds indeed provided $X$ contains a dense open $G$-invariant affine subset.
\end{lem}
\begin{proof} If $G$ is abelian then $G$ is a quotient of the abelian affine algebraic group $U_1\times\ldots\times U_k$, 
which implies the assertion in this case. Otherwise there is a pair of non-commuting $\GG_{\mathrm{a}}$-subgroups $U_i$ and $U_j$. If the general orbits of $U_i$ and $U_j$ have distinct closures then for  the general point $P\in X$ the rational map
\[\A^2_\kk\dasharrow X,\quad (s,t)\mapsto U_i(s)U_j(t)(P)\]
is dominant. So $X$ is unirational. By the Castelnuovo rationality theorem $X$ is rational, 
contrary to our assumption. 

Therefore, all the $U_i$ have the same one-dimensional closures $\bar{\mathscr{O}}$ of general orbits, and these are the ones of $G$. Let $\tilde{\mathscr{O}}$ be the normalization of $\bar{\mathscr{O}}$.
Then $U_i|_{\tilde{\mathscr{O}}}$ and $U_j|_{\tilde{\mathscr{O}}}$ are non-commuting $\GG_{\mathrm a}$-subgroups of $\Aut(\tilde{\mathscr{O}})$, cf.~\cite[Exemple~3]{Dem70}.  
By Lemma~\ref{lem:two actions} $\tilde{\mathscr{O}}\cong\PP^1_\kk$, the fixed points of $U_i|_{\tilde{\mathscr{O}}}$ and $U_j|_{\tilde{\mathscr{O}}}$ are distinct and 
$\langle U_i|_{\tilde{\mathscr{O}}},\,U_j|_{\tilde{\mathscr{O}}}\rangle\cong {\rm PGL}(2,\kk)$. 
Hence $\langle U_i,\,U_j\rangle$ surjects onto ${\rm PGL}(2,\kk)$.

Assume now that $X$ contains a dense open $G$-invariant affine subset $Y$. Then each $U_i$ is transitive  on the general orbit $\mathscr{O}\cong\AA^1_\kk$ of $G$ and the restrictions $U_i|_{\mathscr{O}}$ and $U_j|_{\mathscr{O}}$ commute. Hence also $U_i$ and $U_j$ commute for any $i,j$. 
\end{proof}
\begin{proof}[Proof of Proposition~\ref{thm:proj-surf}] 
By Theorem~\ref{thm:CU}  $\Bir(X)$ verifies the Tits alternative. 
Since the connected subgroup $G\subset\Bir(X)$ contains no nonabelian free subgroup it is solvable. 
If $X$ is irrational then
by Lemma~\ref{lem:irrat} case (1) of the proposition holds. 
Assume now that $X$ is a rational surface and $G\subset\Bir(X)\cong\Bir(\PP^2_\CC)$. By Theorem~\ref{thm:DU-solv} 
$G$ is one of the groups in (1)-(4) of this theorem. Any $\GG_{\mathrm{a}}$-subgroup $U_i$ is infinitely divisible. 
Hence it admits no nontrivial homomorphism to 
the group $\Gamma=H\rtimes \ZZ$ where $H$ is a  finite group. The same holds as well for $G=\langle U_1,\ldots,U_k\rangle$. 
Therefore $G$ cannot be of type (4) in Theorem~\ref{thm:DU-solv}. It cannot be of type (3) either, 
since otherwise the image of $G$ in $N(T)$ is contained in $T=N(T)^0$,  which is impossible. If $G$ is of type (2), that is, $G$ is conjugate to a subgroup of an affine algebraic group, then by Lemma~\ref{lem:alg} $G$ is a unipotent affine algebraic group. 
By Theorem~\ref{thm:D-nilp} $G$ is metabelian, which corresponds to case (1) of the proposition. It remains to treat the  possibility for $G$ to be of type (1) in Theorem~\ref{thm:DU-solv}. 

Thus, we assume in the sequel that $G$ is contained in the de Jonqu\`ere subgroup $J={\rm PGL}(2, \CC(t)) \rtimes {\rm PGL}(2, \CC)$ of $\Bir(\PP^2_\CC)$. 
The latter is equivalent to the fact that $G$ preserves a pencil of rational curves on $\PP^2_\CC$, see e.g.~\cite[p.~4]{Des21}.
Let $\varphi\colon J\to  {\rm PGL}(2, \CC)$ be the quotient morphism. The connected solvable subgroup $\varphi(G)$ is contained in the Borel subgroup $B$ of $ {\rm PGL}(2, \CC)$. 
For any $i=1,\ldots,k$ the image $\varphi(U_i)$ 
is contained in the unipotent radical $R_{\mathrm u}(B)\cong\GG_{\mathrm a}$. 
Thus, we have the exact sequence 
\begin{equation}\label{eq:short} 1\to N
=\ker(\varphi)\to G \stackrel{\varphi}{\longrightarrow} \GG_{\mathrm a}=\GG_{\mathrm a}(\CC).
\end{equation}

Assume first that $N=G$, that is,  
$G\subset {\rm PGL}(2, \CC(t))$. By Lemma~\ref{lem:P1xP1} $G$ is abelian.
Hence $G$ is a quotient of $\prod_{i=1}^k U_i$ 
and once again case (1) of the proposition holds.

Let now $N\neq G$. Then $\varphi$ is surjective,  
$G/N\cong\GG_{\mathrm a}$ and we are in case (2) of the proposition.
\end{proof}
The following example shows that the group $G$ in Proposition~\ref{thm:proj-surf} 
can be an infinite-dimensional (hence non-algebraic) metabelian group. 
\begin{exa}\label{exa:inf-dim}
Consider the $\GG_{\mathrm{a}}$-subgroups $U_1$ and $U_2$ of 
$\Bir(\PP^1_\kk\times\PP^1_\kk)\cong \Bir(\PP^2_\kk)$ given in affine coordinates
$(x,y)\in\A^2_\kk\subset \PP^1_\kk\times\PP^1_\kk$ by formulas
\[U_1=\{(x,y)\mapsto (x+\tau, y)\},\quad U_2=\{(x,y)\mapsto (x, y+\frac{\tau }{x})\}, \quad \tau \in\kk.\]
They are generated by the locally nilpotent derivations
\[\p_1=\p/\p x\in\Der(\kk[x,y])\quad\text{resp.}\quad\p_2=x^{-1}\p/\p y\in\Der(\kk(x)[y]).\]
The Lie subalgebra $L\subset\Der(\kk(x)[y])$ generated by $\p_1$ and $\p_2$ is two-step  solvable. In fact,
\[L:=\langle \p_1, \p_2\rangle=\langle L_1,L_2\rangle \quad\text{where}\quad L_1
=\kk\p_1 \quad\text{and}\quad L_2=\{ f(x^{-1})\p/\p y\,|\,f\in u\kk[u]\}\]
are abelian Lie subalgebras.
We have $[L,L_2]\subset L_2$ where $L_2\subset L$ is an abelian ideal.
Furthermore,  $G=\langle U_1,U_2\rangle$ is an infinite-dimensional 
 subgroup of the metabelian group $R_{\mathrm{u}}(\tilde B)  \rtimes R_{\mathrm{u}}(B)\cong \GG_{\mathrm{a}}(\kk(t))\rtimes\GG_{\mathrm{a}}(\kk)$, that is, of the semidirect product of unipotent radicals 
of the Borel subgroups $\tilde B\subset {\rm PGL}(2, \kk(x))$ and $B\subset {\rm PGL}(2, \kk)$. Thus, $G$ is an infinite-dimensional metabelian subgroup of $\Bir(\PP^2_\kk)$. 
\end{exa}
In the next example we exhibit a $\GG_{\mathrm a}$-subgroup $U$ of the semidirect product of the Borel 
subgroups $\tilde B\rtimes B\subset J$ which is not contained in  the semidirect product of the
unipotent radicals $R_u(\tilde B)\rtimes R_u(B)$. 
\begin{exa}\label{exa:fin-dim}
Consider the  $\GG_{\mathrm a}$-subgroup $U\subset\tilde B\rtimes B\subset J$ acting on $\PP^1_\kk\times\PP^1_\kk$ via birational transformations
\[(x,y)\mapsto\left(x+\tau, \frac{x+\tau}{x}
 y+\frac{\tau}{x}\right),\quad \tau \in \kk.\]
 The  infinitesimal generator of $U$ 
is the derivation 
\[\p=\frac{\p}{\p x}+\frac{y+1}{x}\frac{\p}{\p y}
\quad\text{where}\quad \p^{2}(x)=\p^{2}(y)=0.\]
Clearly, $U\not\subset R_u(\tilde B)\rtimes R_u(B)\cong\GG_{\mathrm a}(\kk(x))\rtimes\GG_{\mathrm a}(\kk)$.
\end{exa}
\begin{rem}\label{rem:FH}
It is known that any $\GG_{\mathrm{a}}$-subgroup of $\Bir(\PP^2_\kk)$ is conjugate to the subgroup of translations
$(x,y)\mapsto (x+\tau,y)$, $\tau \in\kk$, see~\cite[Corollary~3.2]{DF91} and \cite[Corollary~5(i)]{Pop17} (for a more general result); cf.\ also \cite[Proposition 85]{FH21}.
\end{rem}
\section{The case of an affine surface} \label{sec:aff-surf}
 The next result is our main theorem.
\begin{thm}\label{thm:aff-surf}
Let $X$ be an affine surface over $\CC$ and
$G$ be a subgroup of 
$\Aut(X)$ generated by a finite collection of $\GG_{\mathrm{a}}$-subgroups $U_1,\ldots,U_k$. 
Assume that $G$ contains no nonabelian free subgroup. Then $G$ is a metabelian unipotent affine algebraic group. 
\end{thm}
\begin{proof} 
By Proposition~\ref{thm:proj-surf} the assertion holds unless $X$ is rational, $G$ preserves 
a pencil $\Lambda$ of rational curves on $X$ and has a solvable normal subgroup $N$ such that $G/N\cong\GG_{\mathrm{a}}$. 
In this remaining case we consider a functorial resolution of singularities 
$\tilde X\to X$. It is $G$-equivariant and the lift to $\tilde X$ of any regular action of an algebraic group on $X$ 
is regular, see~\cite[\S3.4, Proposition~3.9.1 and Theorem~3.36]{Kol07}.
Since $X$ is a surface, one can take for such a resolution the normalization followed by the minimal resolution of singularities. 
Indeed, both of them are functorial. 

Thus, we may and will assume that $X$ is a smooth quasi-projective surface 
and $\Lambda$ is a $G$-invariant pencil of non-complete rational curves on $X$. The fixed components of $\Lambda$ are $G$-invariant. 
Deleting these components we may suppose that  the base locus ${\rm Bs}(\Lambda)$  is a finite set of points. 
These points are fixed by the action of $G$. We can resolve the base points of  $\Lambda$ via a minimal sequence of successive blowups with zero-dimensional centers.
The $G$-action on $X$ lifts to the resolution and the lift of any $U_i$-action remains regular. 
Therefore we may assume that $\Lambda$ is base-points free. 

By our assumption, $\varphi\colon G\to  \GG_{\mathrm{a}}$ in~\eqref{eq:short} is surjective. 
Hence $\varphi|_{U_l}$ surjects onto $\GG_{\mathrm{a}}$ for some $l\in\{1,\ldots,k\}$.  
The effective $\GG_{\mathrm{a}}$-action on the base $\PP^1_\CC$ of $\Lambda$ has a unique fixed point. Let this point be the point $\infty\in\PP^1$.  
If the corresponding member $\gamma_\infty$ of the pencil $\Lambda$ is a nonzero effective divisor then it is $G$-invariant. Replacing $X$ 
by $X\setminus{\rm supp}(\gamma_\infty)$ we may suppose that $\Lambda$ is parameterized by $\A^1_\CC=\PP^1_\CC\setminus\{\infty\}$. 
The base point free pencil  $\Lambda$ defines a morphism $\pi\colon X\to\A^1_\CC$  which is $G$-equivariant with respect to the $G$-action on $X$ and the $\varphi(G)$-action on $\A^1_\CC$. 
The $\varphi(U_{l})$-action on $\A^1_\CC$ being free the $U_l$-action on $X$ is as well  free. The general member $\gamma$ of $\Lambda$ is a smooth affine curve. 
The free $U_l$-action on $X$ yields 
a trivialization of the fibration $\pi\colon X\to\A^1_\CC$ via the $U_l$-equivariant isomorphism
\[\A^1_\CC\times\gamma\stackrel{\cong}{\longrightarrow} X,\quad (\tau,P)\mapsto U_l(\tau)(P),\quad \tau \in\CC,\] where $U_l$ acts on $\A^1_\CC\times\gamma$ 
via translations on the first factor. 
If the affine curve $\gamma$ is not isomorphic to the affine line then any morphism $\A^1_\CC\to\gamma$ is constant. 
Hence there is a unique $\A^1$-fibration on $X$, namely, the one coming from ${\rm pr}_2\colon \A^1_\CC\times\gamma\to\gamma$  
via the above isomorphism $X\cong\A^1_\CC\times\gamma$. 
It follows that  for any $i=1,\ldots,k$ the general orbits of the  $U_i$-action on $X$ are the general fibers of the latter $\A^1_\CC$-fibration. 
In this case $U_i$ and $U_j$ commute for any $i,j$, hence $G$ is abelian and a quotient of $\prod_{i=1}^k U_i\cong\GG_{\mathrm{a}}^k$. 

Suppose $G$ is not abelian. Then $\gamma \cong\A^1_\CC$ and so $X\cong\A^2_\CC$. Now Proposition~\ref{thm:A2} applies and gives the result.
\end{proof}
\begin{rems} 1.
Alternatively, one can use a detailed description of the subgroup of automorphisms of an 
$\A^1$-ruled affine surface which preserve the ruling, see~\cite[Theorems~8.13 and~8.25]{KPZ17}. 

2. As a curiosity, let us observe that any group $G$  generated by a finite collection of $\GG_{\mathrm{a}}$-subgroups  and acting on an affine surface $X$ over $\CC$ with a non-closed orbit must contain a nonabelian free subgroup. Indeed, otherwise by Theorem~\ref{thm:aff-surf} $G$ is unipotent and then any $G$-orbit is closed in $X$ according to a theorem of Rosenlicht. As a simple example consider the  group $G=\SL_2(\CC)$ acting on the affine plane. It is generated by the $\GG_{\mathrm{a}}$-subgroups
\[U_1=\{(x,y)\mapsto (x+sy,y)\}\quad\text{and}\quad U_2=\{(x,y)\mapsto (x,y+tx)\},\quad s,t\in\CC.\]  The origin of $\A^2_\CC$ is a unique closed $G$-orbit and its complement is an open orbit of $G$. 

However, the above condition is not necessary in order that $G$ contains a nonabelian free subgroup. Indeed, adding a third $\GG_{\mathrm{a}}$-subgroup of shifts $U_3=\{(x,y)\mapsto (x+r,y),\,\,\, r\in\CC\}$ one gets a group $\langle U_1, U_2, U_3\rangle$ acting transitively on the plane and  containing  a nonabelian free subgroup.
\end{rems}
\section{Unipotent group actions and  Gizatullin surfaces}\label{sec:unip-grps}
\subsection{Unipotent group actions on affine surfaces}
\begin{prop}\label{prop:unip}  Let $X$ be an affine surface over $\kk$ and $U$ be a unipotent algebraic subgroup of $\Aut(X)$. 
Then one of the following holds: 
\begin{itemize}
\item[(i)] $U$ has an open orbit in $X$. In this case $X\cong\A^2_\kk$, $U$ is transitive on $\A^2_\kk$ and is conjugate to 
a subgroup of the unipotent de Jonqui\`eres subgroup $\Jonq(\A^2_\kk)$. 
\item[(ii)] The general orbits of $U$ on $X$ are one-dimensional. In this case $U$ is abelian.
\end{itemize}
For any natural number $n\ge 2$  there exists a unipotent subgroup $U\subset\Aut(X)$ of dimension 
$\dim U=n$ satisfying (i), respectively, (ii). 
\end{prop}
\begin{proof}
In case (i) the open orbit of the unipotent algebraic group $U$ is closed in $X$, hence coincides with $X$ 
and is isomorphic to the affine plane $\A^2_\kk$.
By Serre's theorem~\cite[Theorem~8]{Ser80} the algebraic group $U$ is conjugate to a subgroup of either 
$\Aff(\A^2_\kk)$ or $\Jonq(\A^2_\kk)$. In any case,
being unipotent $U$ is conjugate to a subgroup of $\Jonq(\A^2_\kk)$. 
It remains to notice that the unipotent algebraic group $G=\langle U_1,U_2\rangle$ 
in Remark~\ref{rem:der-length} has dimension $n+2$ and acts transitively on the plane $\A^2_\kk$. 

In case (ii) consider a $\GG_{\mathrm{a}}$-subgroup $U_0$ of $U$, 
and let $\p_0\in\lie(U_0)$ be the locally nilpotent derivation generating $U_0$.  
The Lie algebra $\lie(U)$ consists of regular vector fields on $X$ tangent to the orbits of $U$. 
Any such derivation has the form $f\p_0$ where $f$ is a $U$-invariant rational function on $X$. 
It follows that $\lie(U)$, and as consequence also $U$, are abelian.

The last assertion in this case can be confirmed as follows. 
Starting with a $\GG_{\mathrm{a}}$-subgroup $U_0$ of $\Aut(X)$ and a nonconstant $U_0$-invariant function 
$f\in\mathcal{O}(X)$ 
and adding the $\GG_{\mathrm{a}}$-subgroups $U_i=\exp(f^i\p_0)$, $i=1,\ldots,n-1$ 
one obtains an $n$-dimensional abelian unipotent subgroup $U\subset\Aut(X)$ where
\[U=\prod_{i=0}^{n-1} U_i=\{\exp(p(f)\p_0)\,|\,p\in\kk[t],\,\,\,\deg(p)\le n-1\}.\]~\end{proof}
\subsection{Gizatullin surfaces}
We recall a classification of normal affine surfaces with respect to the `size' of the collection of 
$\GG_{\mathrm{a}}$-subgroups of their automorphism groups, see e.g.~\cite[Section~5]{KPZ17}. 
Given such a surface $X$ over $\kk$ we let $\SAut(X)$ be the subgroup of $\Aut(X)$ generated by 
all the $\GG_{\mathrm{a}}$-subgroups. 
The \emph{Makar-Limanov invariant} $\ML(X)$ measures the complexity of the action of the group $\SAut(X)$ on $X$. 
More precisely, we say that $X$ is of class $\ML_i$ where $i\in\{0,1,2\}$ if the codimension in $X$ 
of the general $\SAut(X)$-orbit is equal to~$i$. 
The surface $X$ is said to be \emph{rigid} if  it is of class $\ML_2$, that is, does not admit any nontrivial 
$\GG_{\mathrm{a}}$-action. 
If $X$ is of class $\ML_1$ then all the effective $\GG_{\mathrm{a}}$-actions on $X$ have the same general orbits 
contained in the fibers of an $\A^1_\kk$-fibration $X\to B$ over a smooth affine curve $B$. 
In this case, under a mild addition assumption, $\SAut(X)$ is an abelian nested ind-group~\cite[Section~6]{KPZ17}. 
A simple example is provided by the surface $X=(\A^1_\kk\setminus\{0\})\times \A^1_\kk$. Indeed, this surfaces carries the unique $\A^1_\kk$-fibration ${\rm pr}_1\colon X\to \A^1_\kk\setminus\{0\}$. 

The most interesting class is the class $\ML_0$. It consists of the Gizatullin surfaces. 
A  normal affine surface $X$ over $\kk$ non-isomorphic to $(\A^1_\kk\setminus\{0\})\times\A^1_\kk$ is called 
a \emph{Gizatullin surface} if $X$
can be embedded into a complete surface $\bar X$ by adding a reduced boundary divisor $D=\bar X\setminus X$ which is a chain of smooth rational curves with simple normal crossings. By the Gizatullin theorem~\cite[Theorems~2 and~3]{Giz71} a normal affine surface $X$  is Gizatullin if and only if $\Aut(X)$ contains two  $\GG_{\mathrm{a}}$-subgroups $U_1,U_2$ whose general orbits are different, in other words, if the group $\SAut(X)$ acts on $X$ with an open orbit, that is, $\ML(X)=\kk$.
Given a Gizatullin surface $X$ the group $\SAut(X)$ acts infinitely transitively on its open orbit~\cite[Example~2.3]{AFKKZ13}. For any natural number $n>1$ there exists a smooth Gizatullin surface $X$ carrying an $n$-parameter family of $\A^1_\kk$-fibrations over $\A^1_\kk$ such that any two of them are not equivalent  under the $\Aut(X)$-action on $X$~\cite[Example~6.3.21]{FKZ11}. Such a surface $X$ also carries an $n$-parameter family of $\GG_{\mathrm{a}}$-subgroups acting along these $\A^1_\kk$-fibrations and pairwise non-conjugate in $\Aut(X)$. 
Notice that for  $X=\A^2_\kk$ any two $\A^1_\kk$-fibrations over $\A^1_\kk$ are equivalent up to the $\Aut(\A^2_\kk)$-action on $\A^2_\kk$; this follows e.g. from Rentchler's Theorem~\ref{thm:R}.

Using Theorem~\ref{thm:main} and Proposition~\ref{prop:unip} we deduce the following result. 
\begin{prop}\label{new-cor} Let $X$ be a Gizatullin surface over $\CC$ non-isomorphic to $\A^2_\CC$ and let $G$ be a subgroup of $\Aut(X)$ generated by a collection of $\GG_{\mathrm{a}}$-subgroups. 
Then the following Tits-type alternative holds: either $G$ is abelian and acts on $X$ with one-dimensional general orbits, 
or $G$ contains a nonabelian free subgroup and acts with an open orbit on $X$. 
\end{prop}
\begin{proof} If  $G$ acts on $X$ with one-dimensional general orbits than it preserves fiberwise an $\A^1_\CC$-fibration 
$X\to\A^1_\CC$. Hence any two $\GG_{\mathrm{a}}$-subgroups 
of $G$ commute and so $G$ is abelian. 

Otherwise $G$ contains two non-commuting $\GG_{\mathrm{a}}$-subgroups
$U_1$ and $U_2$ which act on $X$ with different general orbits. The subgroup 
$U=\langle U_1, U_2\rangle$ of $G$ acts on $X$ with an open orbit. If $U$ is unipotent then by
Proposition~\ref{prop:unip} 
$X\cong\A^2_\CC$ contrary to our assumption. Therefore, by Theorem~\ref{thm:main} $U$ contains a nonabelian free subgroup.
\end{proof}
\section{Note added in proof}\label{sec:note}
We provide here an alternative (though less elementary) proof of Theorem \ref{thm:aff-surf} which runs over more general fields.  
\begin{thm}\label{thm:aff-surf-kk} Fix an algebraically closed field $\kk$ of characteristic zero. Assume that the Tits alternative for finitely generated subgroups holds  in $\Bir(\PP^2_\kk)$.
Let $X$ be an affine surface over $\kk$ and $G$ be a subgroup of 
$\Aut(X)$ generated by a finite collection of $\GG_{\mathrm{a}}$-subgroups $U_1,\ldots,U_k$. 
If $G$ contains no nonabelian free subgroup then $G$ is a metabelian unipotent affine algebraic group. 
\end{thm}
 In the proof of Theorem \ref{thm:aff-surf-kk} we use the following recent result, see \cite[Theorem B]{KZ22}. 
\begin{thm}\label{thm:KZ} Let $X$ be an affine algebraic variety defined over an algebraically closed field $\kk$ of characteristic zero.
Let $G \subset \Aut(X)$ be a  subgroup generated by a finite family $U_1,\ldots,U_n$ of connected algebraic subgroups. If $G$ is solvable, then $G$ is an algebraic group. If in addition the $U_i$ are unipotent, then $G$ is a unipotent algebraic group.
\end{thm}
\begin{proof}[Proof of Theorem \ref{thm:aff-surf-kk}] 
If $X$ is irrational then by Lemma \ref{lem:irrat} $G$ is an abelian unipotent algebraic group. Suppose further that $X$ is rational and that $G$ contains no nonabelian free subgroup. 
For any nonunit element $u_i\in U_i$ the subgroup $\langle u_i\rangle$ is Zariski dense in $U_i$. 
Hence the finitely generated subgroup $H=\langle u_1,\ldots,u_n\rangle$ of $G$ is dense in $G$ for the Zariski ind-topology on $\Aut(X)$, see \cite[Definition 1.1.3 and Theorem 5.1.1]{FK18}. Since $H$ contains no nonabelian free subgroup, by our assumption it is virtually solvable. Let $H_0\subset H$ be a solvable subgroup of finite index $m$. Then $u_i^{m_i}\in H_0$ for some $m_i\in\{1,\ldots,m\}$. It follows that $H_0\cap U_i$ is Zariski dense in $U_i$, and so $H_0$ is Zariski dense in $G$. We conclude that $G=\overline{H_0}$ is solvable. 

Notice that $X$ and the $U_i$ are defined over a subfield $\kk_0\subset\kk$ which is a finite extension of $\QQ$. Fix an embedding $\kk_0\hookrightarrow\CC$ and consider the complex models $X_\CC$ of $X$ and $U_i(\CC)\cong\GG_a(\CC)$ of the $U_i$ where $U_i(\CC)$ acts regularly on $X_\CC$. Consider also the subgroup $G(\CC)=\langle U_1(\CC),\ldots, U_n(\CC)\rangle$ of $\Aut(X_\CC)$.
By Theorem \ref{thm:KZ} $G$ is a unipotent algebraic group.  
Hence $G(\CC)$ is a unipotent algebraic group acting on $X_\CC$. Since $G(\CC)\subset \Bir(\PP^2_\CC)$ is a nilpotent subgroup, by Theorem \ref{thm:D-nilp} $G(\CC)$ is a metabelian  connected algebraic group. Then $G\subset G(\CC)$ is also metabelian. 
\end{proof}

The same arguments prove the following multidimensional analog of Theorem \ref{thm:aff-surf-kk}.
\begin{thm} \label{thm:AZ-KZ} Let $X$ be an affine algebraic variety defined over an algebraically closed field $\kk$ of characteristic zero, and let a subgroup $H\subset \Aut(X)$ verify the Tits alternative for finitely generated subgroups. 
Let $G \subset H$ be a  subgroup generated by a finite family $U_1,\ldots,U_n$ of unipotent algebraic subgroups. If $G$ contains no nonabelian free subgroup then $G$ is a unipotent algebraic group. 
\end{thm}
Using this theorem and \cite[Main theorem]{LP22} resp. \cite[Theorem C]{BFL14} we come to the following 
\begin{cor}\label{cor-LP} 
The conclusion of Theorem \ref{thm:AZ-KZ} holds for 
\begin{itemize}\item
$X=\A^3_\kk$ and $H={\rm Tame}(\A^3_\kk)$ being the tame subgroup of $\Aut(X)$;
\item $X=\SL(2,\CC)$ viewed as an affine threefold and $H={\rm Tame}(\SL(2,\CC))$ being the tame subgroup of $\Aut(X)$. 
\end{itemize}
\end{cor}
\begin{rem}\label{rem:Lamy-book} According to \cite[Theorem 16.2]{Lam22} the group $\Bir(\PP^2_k)$ satisfies the Tits alternative for any field $k$ of characteristic zero. In particular, the conclusion of Theorem \ref{thm:aff-surf-kk} holds  over any algebraically closed field of characteristic zero. Besides, this gives a hope that also Theorem \ref{thm:DU-solv} and Proposition \ref{thm:proj-surf} hold in this generality.  

The positive characteristic case of Conjecture \ref{conj} was added by a suggestion of the anonymous referee. A later discussion with St\'ephane Lamy gives hope that the algebraic closedness assumption of the base field is also superfluous.
This is supported by the fact that the Tits alternative for finitely generated subgroups holds for $\Aut(\AA^2_k)$ over any prime characteristic field $k$, see \cite[Proposition 16.16]{Lam22}. Moreover, an analog of Theorem \ref{thm:L} also holds for $\Aut(\AA^2_k)$ over any field $k$ according to \cite[Proposition 7.49]{Lam22}.

Theorem 16.25 in \cite{Lam22} claims that the Tits alternative for finitely generated subgroups 
holds in the Cremona group $\Bir(\PP^2_k)$ over a prime characteristic field $k$. The proof of this theorem in \cite{Lam22} is under construction, with some arguments just sketched. See also the announcement \cite[Theorem 6.1]{Can18} with a reference to \cite{Can11}; however, \cite{Can11} deals exceptionally with complex surfaces and uses transcendental tools. Notice that there exists
a classification of maximal algebraic subgroups of $\Bir(\PP^2)$ over a perfect field; see  \cite{SZ21}.

Our proof of Theorem \ref{thm:aff-surf} uses the Hironaka resolution of singularities. It is known 
to hold for excellent schemes of dimension 2 defined over a prime characteristic field, see e.g. \cite{Lip78}. Besides, we exploit the expression of a $\GG_{\mathrm{a}}$-subgroup via a locally nilpotent derivation, see e.g. Example \ref{exa:inf-dim} 
and the proof of Proposition \ref{prop:unip}. Such an expression is not available in positive characteristic. 
\end{rem}

\end{document}